\documentclass[letterpaper,10pt]{amsart}

\usepackage{amsfonts}
\usepackage{amsmath}
\usepackage{amssymb}
\usepackage{amsthm}
\usepackage{subfigure}
\usepackage[T1]{fontenc}
\usepackage[dvips]{psfrag,graphics}
\usepackage{enumerate}

\newtheorem{theorem}{Theorem}[section]

\newtheorem{proposition}[theorem]{Proposition}
\newtheorem{definition}{Definition}
\newtheorem{lemma}[theorem]{Lemma}
\newtheorem{corollary}[theorem]{Corollary}

\newenvironment{listi}{\begin{enumerate}[\upshape(i)]\setlength{\itemsep}{1pt}\setlength{\parskip}{0pt}\setlength{\parsep}{0pt}}
{\end{enumerate}}

\newenvironment{lista}{\begin{enumerate}[\upshape(a)]\setlength{\itemsep}{1pt}\setlength{\parskip}{0pt}\setlength{\parsep}{0pt}}
{\end{enumerate}}

\newcommand{\Ha}{\mathcal{H}}
\newcommand{\la}{\lambda}

\newcommand{\al}{\alpha}

\newcommand{\ep}{\epsilon}

\newcommand{\om}{\omega}
\newcommand{\Om}{\Omega}

\newcommand{\diam}{\text{diam}\,}

\DeclareMathOperator{\fav}{Fav}

\begin{document}
\title{On planar self-similar sets with a dense set of rotations}

\author{Kemal Ilgar Ero\u{g}lu}
\address{Department of Mathematics\\
	University of Washington\\
	Box 354350\\
	Seattle, WA 98195-4350}
\email{kieroglu@math.washington.edu}


\date{March 8, 2006}
\maketitle

\begin{abstract}
We prove that if $E$ is a planar self-similar set with similarity dimension $d$ whose defining maps generate a dense set of rotations, then the $d$-dimensional Hausdorff measure of the orthogonal projection of $E$ onto any line is zero. We also prove that the radial projection of $E$ centered at any point in the plane also has zero $d$-dimensional Hausdorff measure. Then we consider a special subclass of these sets and give an upper bound for the Favard length of $E(\rho)$ where $E(\rho)$ denotes the $\rho$-neighborhood of the set $E$.
\end{abstract}

\section{Introduction}

In this paper we investigate the orthogonal and radial projection properties of some self-similar sets in the plane. Planar self similar sets are the attractors of iterated function systems whose maps are contracting similitudes. By a classical result of Marstrand \cite{Marstrand}, if $E$ is a planar Borel set with $s:=\dim_H E \leq 1$ (where $\dim_H$ denotes the Hausdorff dimension) then in Lebesgue almost all directions the orthogonal projections onto lines have Hausdorff dimension $s$. If $\dim_H E>1$ then almost all projections have positive Lebesgue measure. Therefore the natural question is to ask when the projections have positive $s$-dimensional Hausdorff measure in the case when $s\leq1$.

It is useful to partition planar self-similar sets sets into two categories when studying their orthogonal projections onto lines, namely, the case when the similitudes do not involve rotations or reflections, and the others. The sets whose defining maps do not involve rotations are distinguished from the others with the relatively simple structure of their projections: These projections are self-similar sets in $\mathbb{R}$. The defining maps for the projections can be viewed as a family of maps depending on a parameter (the projection angle) and measure theoretic arguments can be made about the properties of ``typical'' projections. We refer to \cite{SolScan} for the details and some applications of this method.

In the case when rotations are involved, the projections are no longer self-similar, thus it is significantly more complicated to study their structure. Our main result is concerned with the case when the defining maps generate a dense set of rotations, and partly answers a question by Mattila in \cite{MattilaSurvey}. This would be the case when, for example, one of the maps involved rotation by an irrational multiple of $\pi$. In fact, if none of the maps involve reflections, our condition is equivalent to having at least one map with an irrational rotation. The idea of the proof is that, if the set of rotations is dense, there are many groups of smaller copies of the original set that are approximately the same size and aligned in a dense set of directions. This is shown by modifying a ``doubling'' argument that has been used in \cite{SolNonlinear}. 
These copies ``pile up'' above a typical point in the projection, making the density of the projected measure infinite. Density of directions suggests that this is the case in a dense set of directions, and hence, by approximation, all directions. We prove the following:

{\bf Theorem.\ }\emph{Let $E\subset\mathbb{R}^2$ be a self-similar set whose defining maps generate a dense set of rotations modulo $2\pi$. If $\gamma$ is the similarity dimension of $E$, then for all lines $l$, the $\gamma$-dimensional Hausdorff measure of the orthogonal projection of $E$ onto $l$ is zero.}

The case when the rotations are in a discrete set of directions is still an open question.

Finally, in the last section we give an upper bound for the Favard length of the $\rho$-neighborhood $E(\rho)$ of $E$ where $E$ is a homogeneous self-similar set of similarity dimension $1$ whose defining maps include two non-rotating maps and a (non-reflecting) rotation by a Diophantine multiple of $\pi$. A number $\alpha$ is called Diophantine if there exist $c,d>0$ such that $|N\alpha - M|> c N^{-d}$ for any two integers $N$ and $M$ (we say the number is $(c,d)$-Diophantine in this case). Recall also that a self-similar set is called homogeneous if all defining maps have the same contraction factor. The Favard length of a planar set $E$ is given by
\[
\fav (E) = \int_0^\pi \mathcal{L}^1 (E^{\theta}) d\theta
\]
where $\mathcal{L}^1 (E^{\theta})$ is the Lebesgue measure of the orthogonal projection of $E$ onto $l_\theta$, the line through the origin making angle $\theta$ with the positive $x$-axis. Besicovitch's projection theorems tell us that in the above case, the Favard length of $E(\rho)$ converges to zero since $E$ is irregular, but in very few cases we have precise estimates for the decay rate. A set $E$ is called a 1-set if $0<\Ha^1(E)<\infty$. There is the known lower bound $c/(-\log \rho)$ when $E$ is any Borel 1-set in the plane (see \cite{Mattila_Favard}).

We will prove the following theorem:

{\bf Theorem.\ }\emph{Let $E$ be the attractor of a homogeneous iterated function system with similarity dimension 1 which produces two non-rotating maps of the same contraction factor and a map rotating by angle $\theta_1$, where $\theta_1/{2\pi}$ is $(c,d)$-Diophantine. Then, denoting the orthogonal projection of $E(\rho)$ onto $l_\theta$ by $E^{\theta}(\rho)$, for any $\delta>0$ there exists $A,B>0$ such that
\begin{equation}
\mathcal{L}^1 (E^{\theta}(\rho)) \leq \frac{A}{(\log(-\log \rho))^B}
\end{equation}
uniformly for all $\theta$, thus
\[
\fav (E(\rho)) \leq \frac{\pi A}{(\log(-\log \rho))^B}.
\] 
}
(This theorem will be stated more precisely in section \ref{sec_favard} (see Theorem \ref{MainFav}).

In \cite{Sol_Buffon}, it was proven that for self-similar sets of similarity dimension $1$ with strong separation and without rotation, the bound
\begin{equation}\label{log_star}
\fav E(\rho) \leq C \exp \left( -a \log_{*} \left(\frac{1}{\rho}\right) \right)
\end{equation}
holds for some $C,a>0$, where
\begin{equation}\label{def_log_star}
\log_{*}  x := \min \left\lbrace n \geq  0:\ \underbrace{\log\log\ldots\log}_{n} x \leq 1 \right\rbrace.
\end{equation}

\begin{figure}[t]\label{fig_1}

\centering
\subfigure{
\scalebox{.5}{\includegraphics{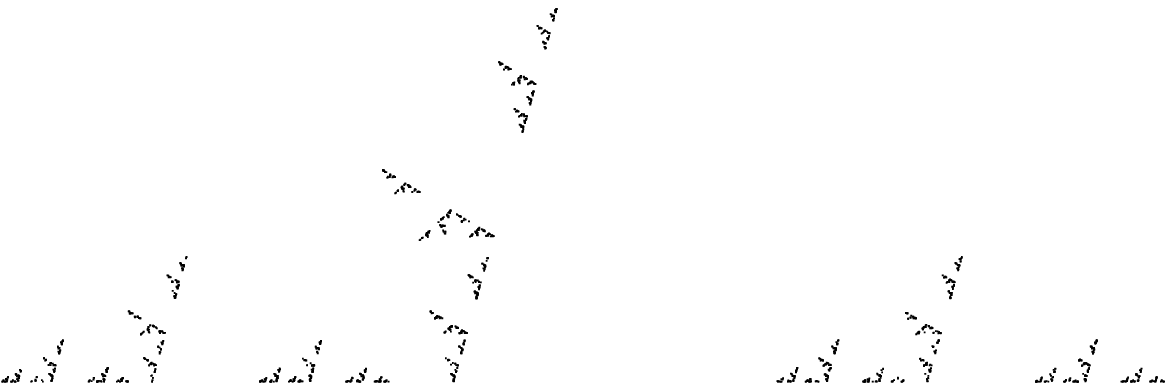}}
}
\qquad \qquad
\subfigure{
\scalebox{.5}{\includegraphics{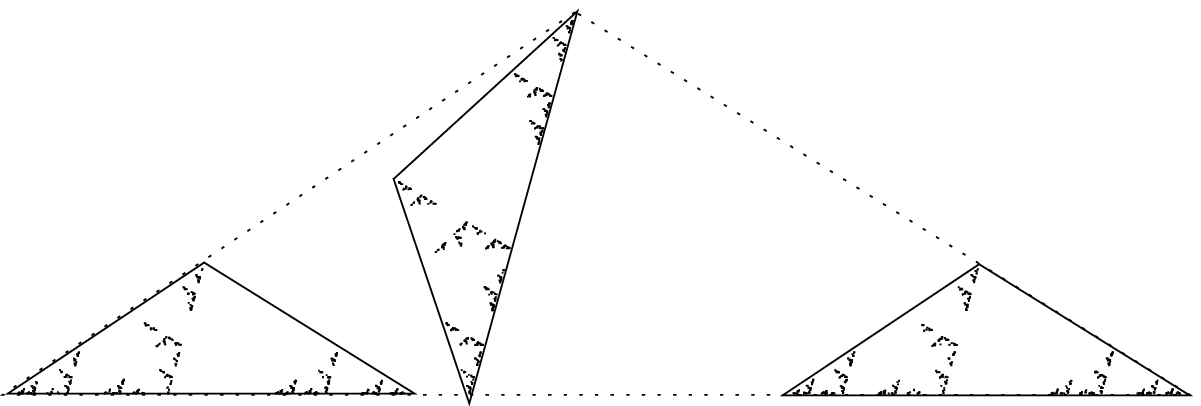}}
}
\caption{}
\end{figure}

The class of sets we study in this paper are the only self-similar sets that are currently known to obey a bound better than (\ref{log_star}).

\medskip
\noindent \emph{Example. } In Figure \ref{fig_1} is shown a homogeneous self-similar set with dimension 1. There are three defining maps two which don't rotate and the third map rotates by $(1+\sqrt{2})\pi$ radians. The number $1+\sqrt{2}$ is Diophantine with $d=2$ since it is a quadratic irrational (see e.g.~\cite{Baker}). For this set we can take $B=\frac{\log2}{3.3 \log 3}$.

\section{Projections of planar self-similar sets with a dense set of rotations}\label{prob2}

We first introduce the terminology before we state the result:

An iterated function system is a finite collection of Lipschitz maps on $\mathbb{R}^2$ with Lipschitz constants less than 1.

Consider an IFS in the plane with maps $F_1, \ldots,F_m$. There is a unique nonempty compact set $E$, called the attractor of this iterated function system, that satisfies
\begin{equation}\label{attractor}
E= \bigcup_{i=1}^m F_i(E).
\end{equation}
We say $E$ is self-similar if the maps $F_i$ are similitudes. The set $E$ can be viewed as the image of a ``projection'' $\Pi$ from a symbol space \mbox{$\Om=\{1,\ldots,m\}^\mathbb{N}$}, given by 
\[
\Pi:\ i_1 i_2 i_3 \ldots \longrightarrow \lim_{n\to \infty} F_{i_1} \circ F_{i_2} \circ \cdots \circ F_{i_n} (E).
\]
The limit is a singleton since the maps are contractions. Finite sequences $u=u_1 u_2 \cdots u_n$ of symbols will be called words. We will use $[u]$ to denote the set of sequences in $\Om$ starting with the word $u$ (such sets are called cylinder sets).

Now assume $E$ is self-similar. Let $l_\theta$ be the line through the origin making an angle of $\theta$ with the positive $x$-axis and let $\Pi_\theta$ be the composition of $\Pi$ with the orthogonal projection onto $l_\theta$. Assume each $F_i$ has contraction rate $r_i$. The similarity dimension of the system above is the number $\gamma$ such that
\[
\sum_{i=1}^m r_i^\gamma=1.
\]
We will shortly call $\gamma$ the similarity dimension of $E$, when the iterated function system in question is clear from the context (Note that $E$ can be produced by different iterated function systems). It is well-known that $\dim_H E \leq \gamma$. For many purposes the interesting situation is the case when $\mathcal{H}^\gamma(E)>0$. The Bandt-Graf condition \cite{BG} and the open set condition (see \cite{Fal2}) are examples of necessary and sufficient conditions for this to be true.

For a word $u=u_1\cdots u_n$ define $r_u=r_{u_1}\cdots r_{u_n}$. Define
\[
F_u := F_{u_1}\circ\cdots\circ F_{u_n}.
\]
We denote by $\bar{u}$ the infinite sequence $uuu\ldots$ The length of $u$ will be denoted by $|u|$.

Each $F_i$ is either a scaling and reflection composed with a translation, or a scaling and rotation composed with translation. Given a word $u$ we can write

\[
F_u \left(\begin{bmatrix} x \\ y\end{bmatrix} \right) = r_u \begin{bmatrix} \cos\theta_u & - O_u \sin\theta_u \\ \sin\theta_u & O_u\cos\theta_u \\ \end{bmatrix} \begin{bmatrix}x \\ y \end{bmatrix} + \begin{bmatrix} x_u \\ y_u \end{bmatrix}
\]
where $O_u=-1$ if $F_u$ contains a reflection (about the line $\theta=\theta_u /2$), or $1$ if it contains a rotation only (by the angle $\theta_u$). Note that if $O_u = O_v =1$ then $\theta_{uv} = \theta_u + \theta_v$, and $|\theta_u -\theta_v|=|\theta_{\al u} -\theta_{\al v}|$ for any words $u,v,\al$. In this paper we will consider angles modulo $2\pi$, e.g.~we will write $|\theta_u -\theta_v|$ for \mbox{$|(\theta_u -\theta_v) \mod2\pi|$}.
We are assuming that the set 
\[
\Theta=\{\theta_u\ :\ \text{$u$ is a word}\}
\]
is dense modulo $2\pi$. Our result is as follows:

\begin{theorem}\label{Main2}
Let $E\subset\mathbb{R}^2$ be a self-similar set such that $\Theta$ is dense modulo $2\pi$. If $\gamma$ is the similarity dimension of $E$, then for all lines $l$, the $\gamma$-dimensional Hausdorff measure of the orthogonal projection of $E$ onto $l$ is zero.
\end{theorem}

Note that if $\Theta$ is dense, then $\{\theta_u\ :\ O_u=1\}$ is also dense. Given an $\ep>0$ we will denote by $a=a(\ep)$ some fixed word such that $O_a=1$ and $|\theta_a|<\ep$.

\bigskip
Let $\gamma$ be the similarity dimension of $E$. We define $\mu$ to be the product measure $\mu=\{r_1^\gamma,...,r_m^\gamma \}^{\mathbb{N}}$ on $\Om$, that is $\mu([u])= r_u^\gamma$. Let $D=\diam E$.

\begin{definition}\label{def_ep}
Two words $u,v$ are called $(\ep,\theta)$-relatively close if
\begin{listi}
\item $\frac{r_u}{r_v} \in (e^{-\ep},e^\ep)$;
\item $|\theta_u - \theta_v| < \ep$ and $O_u=O_v$;
\item There exists $\om=\om(u,v)\in\Om$ such that $|\Pi_\theta (u\om) - \Pi_\theta (v\om)| < \ep D \min\{r_u,r_v\}$.
\end{listi}
\end{definition}

This definition means the smaller copies of the set $E$ defined by $u$ and $v$ have relatively the same size and orientation, and the convex hulls of their projections on $l_\theta$ have large overlap.

The proof of Theorem \ref{Main2} is based on a sequence of lemmas:

\begin{lemma} \label{cl1of2}
Given any $\ep>0$, an angle $\phi$ and word $u$, for $\mu$-almost all points $\om \in \Om$ there exist infinitely many words $s_j$ such that $s_j u$ is a prefix of $\om$, $O_{s_j}=1$ and $|\theta_{s_j} - \phi| < \ep$.
\end{lemma}
\begin{proof}
Let $a=a(\ep)$. There is $N\in\mathbb{N}$ such that the set

\[
\{\theta_a,\theta_{aa},\theta_{aaa},\ldots,\theta_{a^N}\}
\]
forms an $\ep$-net modulo $2\pi$. Given any word $v$, if $O_v=1$ we can concatenate a sequence $t_v$ of $a$'s of length no more than $N$ so that $O_{vt_v}=1$ and $|\theta_{vt_v} - \phi|<\ep$. If there is an orientation reversing $F_{i_0}$ and $O_v=-1$ then we concatenate to $v{i_0}$ a suffix $t_v$ of $a$'s as above to get $O_{v{i_0}t_v}=1$ and $|\theta_{v{i_0}t_v} - \phi|<\ep$. So we conclude that there is an integer $N_0=N_0(\ep)$ such that given any word $v$ we can find a suffix $t_v$ of length no more than $N_0$ satisfying $O_{vt_v}=1$ and $|\theta_{vt_v} - \phi|<\ep$.

Given this $N_0$, there exists $c>0$ such that if $t$ is a word of length no more than $N_0$ then $\mu([t u]) \geq c>0$.

Let $A$ be the set of points in $\Om$ where the claim fails. Given $n \in \mathbb{N}$, let $A_n$ be the points $\om$ of $A$ for which no prefix $s_j$ of length $\ge n$ exists as in the claim. Then $A = \cup A_n$. Therefore it suffices to prove that $\mu(A_n)=0$ for any $n$.

Now we fix $n$. Write $\Om$ as a disjoint union of cylinders represented by words length bigger than $n$ and let $S_0$ be the set of words corresponding to the cylinders in this union. For each word $v\in S_0$, let $t_v$ be a word of length no more than $N_0$ such that $O_{vt_v}=1$ and $|\theta_{v t_v} - \phi|< \ep$. Let $\Om_1 = \Om \setminus \cup_{v\in S_0} [v t_v u]$. Then $\mu(\Om_1) < 1-c$.
Given $\Om_i$, define $\Om_{i+1}$ in the same way: Write $\Om_i$ as a disjoint union of cylinders and let $S_i$ be the set of corresponding words. For each word $v$ in $S_i$, find a word $t_v$ as above. Define $\Om_{i+1} = \Om_i \setminus \cup_{v \in S_i} [v t_v u]$.
Note that $\mu(\Om_{i+1}) < (1-c) \mu(\Om_i)$ for all $i$. Since $A_n$ is a subset of $\cap \Om_i$ and $\mu(\cap \Om_i)=0$, the result follows.
\end{proof}

The following corollary will be useful when studying visibility properties: 

\begin{corollary}\label{corol_to_1of2}
Let $\omega \mapsto \theta_\omega \in [0,2\pi)$ be any function. Then, given a word $u$ and $\ep>0$, for $\mu$-almost all $\om\in\Om$ there exist prefixes $s_j$ such that $s_j u$ is a prefix of $\om$, $O_{s_j}=1$ and $|\theta_{s_j} - \theta_\om|<\ep$.
\end{corollary}
\begin{proof}
Let $\{\phi_n\}$ be a finite collection of angles that form an $\ep/2$-net modulo $2\pi$. Then, given a fixed $\phi_n$, by Lemma~\ref{cl1of2} $\mu$-almost all $\om\in\Om$ have prefixes $s_j$ such that $s_j u$ is a prefix of $\om$, $O_{s_j}=1$ and $|\theta_{s_j} - \theta_{\phi_n}|<\ep/2$. Since $\{\phi_n\}$ forms an $\ep/2$-net, this implies that $\mu$-almost all $\om\in\Om$ have prefixes satisfying the conditions of our claim.
\end{proof}

A set of real numbers is called \emph{$\sigma$-arithmetic} if all numbers in the set are integer multiples of $\sigma$ and $\sigma$ is the biggest number with this property. We quote a probabilistic lemma (see \cite{Feller}, Vol II, Lemma V.4.2):
\begin{lemma}[Feller]\label{Feller2}
Let $F$ be a distribution in $\mathbb{R}$ concentrated on $[0,\infty)$ but not at the origin, and $\Sigma$ the set formed by the points of increase of $F, F*F, F*F*F,\ldots$. If $F$ is not arithmetic, then $\Sigma$ is asymptotically dense at infinity in the sense that for given $\ep>0$ and $K$ sufficiently large, the interval $(K,K+\ep)$ contains points of $\Sigma$. If $F$ has $\la$-arithmetic support then $\Sigma$ contains all points $n\la$ for $n$ sufficiently large.
\end{lemma}

\begin{lemma}\label{cl2of2}
Given any $\ep>0$ and a function $\theta \mapsto \phi(\theta)$, there exist $u,v,\theta$ such that $u,v$ are distinct and $(\ep,\theta)$-relatively close with $O_u=O_v=1$. Moreover, in addition to the conditions (i)-(iii) of Definition \ref{def_ep}, we can choose $u,v$ and $\theta$ in such a way that: (iv) $|\phi(\theta) - \theta_u|< \ep$, $|\phi(\theta) - \theta_v|< \ep$, and (v) there exists $\om$ such that $\Pi_\theta (u\om) = \Pi_\theta (v\om)$.
\end{lemma}

\begin{proof}
\begin{figure}[h]
\psfrag{label01}{$l_\theta$}
\psfrag{label02}{$u$}
\psfrag{label03}{$v$}
\psfrag{label04}{$ua$}
\psfrag{label05}{$va$}
\psfrag{labeln1}{$\Pi(u\bar{a})$}
\psfrag{labeln2}{$\Pi(v\bar{a})$}
\scalebox{.8}{\includegraphics{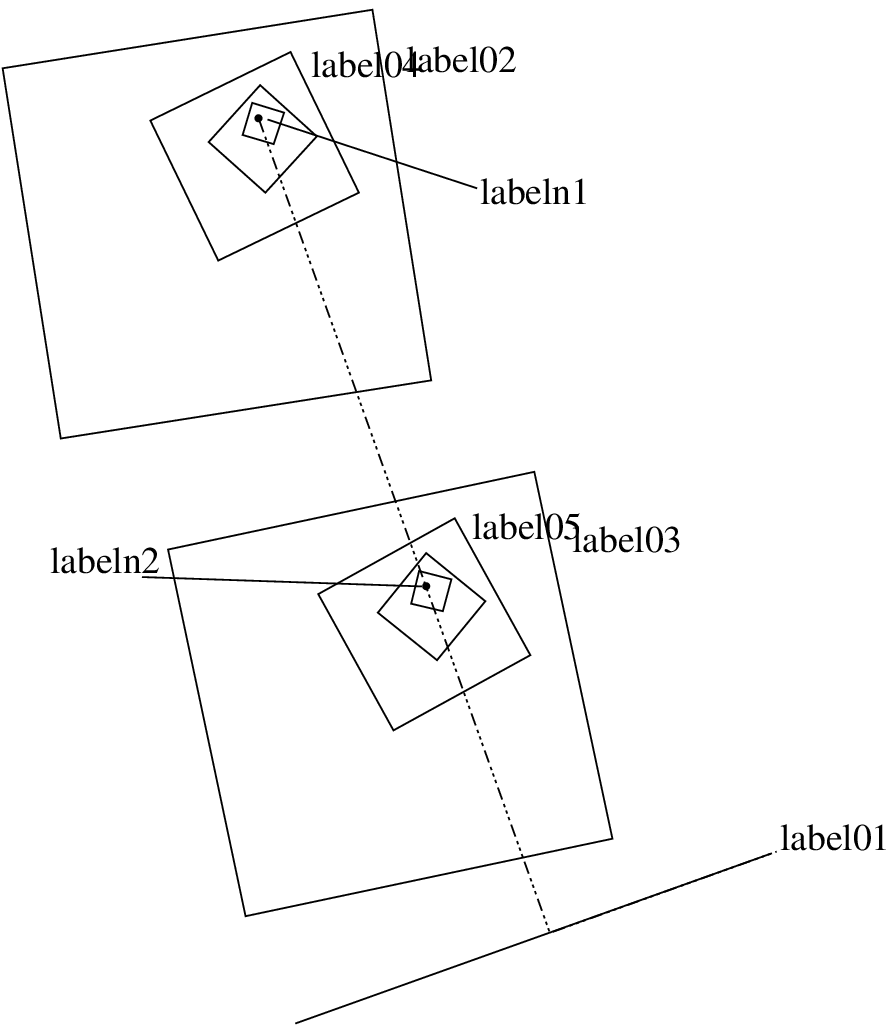}}
\caption{}
\label{resim0}
\end{figure}

Note that (v) trivially implies (iii), so it suffices to check condition (v) alone to prove (iii).
Let $r_{min}=\min_i r_i$. Given $r>0$, let $C_r$ be the the set of all cylinders $[s]$ with 
\[
 r r_{min} < r_s \leq r.
\]
Let $N_r = \# C_r$. We first consider the case when $S=\{-\log r_1 ,\ldots, -\log r_m\}$ is not arithmetic: Let $F$ be a distribution supported on $S$. Find $K$ sufficiently large as in Lemma \ref{Feller2} for $\ep/2$. Let $\tilde{r}$ be such that
\[
-\log \tilde{r} = -\log r -\log r_{min} + K.
\]
Then each $[s]\in C_r$ has a subcylinder $[u]$ (i.e. $s$ is a prefix to $u$) such that $\frac{r_u}{\tilde{r}} \in (e^{-\ep/2}, e^{\ep/2})$. So there exist at least $N_r$ cylinders $u_i$ such that $\frac{r_{u_i}}{r_{u_j}} \in (e^{-\ep}, e^{\ep})$ for each pair $u_i,u_j$. Since $N_r \to \infty$ as $r\to 0$, by choosing $r$ small enough we can also find a pair $\tilde{u},\tilde{v}$ such that $O_{\tilde{u}} = O_{\tilde{v}}$ and $|\theta_{\tilde{u}} - \theta_{\tilde{v}}|<\ep$. If $O_{\tilde{u}}=O_{\tilde{v}}=1$ then let $u=\tilde{u}$ and $v=\tilde{v}$. If $O_{\tilde{u}}=O_{\tilde{v}}=-1$ and $F_{i_0}$ is an orientation reversing map in the iterated function system, let $u=\tilde{u} i_0$ and $v=\tilde{v} i_0$.
Then $O_u=O_v=1$ and $|\theta_u - \theta_v|<\ep$, thus (i) and (ii) in Definition~\ref{def_ep} are satisfied.

In the case when $S$ is $\sigma$-arithmetic for some $\sigma>0$, we again consider a distribution $F$ supported on $S$, and using the lemma find $K$ such that $n\sigma \in \Sigma$ for all $n\geq K$. Let $\tilde{r}$ be given by 
\[
\log {\tilde{r}} = \left[ \frac{-\log r - \log r_{min}}{\sigma} + K +1 \right] \sigma.
\]
Then there exists a subcylinder $[u]$ of $[s]$ with $r_u =\tilde{r}$. This is true for each $[s]\in C_r$. The rest of the argument is as in the non-arithmetic case; we can find $u,v$ satisfying the first two conditions of Definition \ref{def_ep}.

We set $a=a(\ep/2)$. Choosing $\theta$ to be $\pi/2$ plus the direction of the line segment joining $\Pi(u\bar{a})$ to $\Pi(v\bar{a})$ and setting $\om=\bar{a}$, (v) is satisfied (if these two points coincide then any direction can be chosen as $\theta$). For (iv), consider the pairs
\begin{equation*}
(u,v), (ua,va), (uaa,vaa), (uaaa,vaaa),\ldots
\end{equation*}

All these pairs satisfy (i)-(ii) and (v) (hence (iii)) with the same $\ep$, $\theta$ and $\om$. Since $|\theta_a|<\ep/2$, some of them will also satisfy (iv) (see Fig. \ref{resim0}).
 \end{proof}

\begin{lemma}\label{cl3of2}
Given any $N$ and $\ep>0$, there exist distinct words $u_1,\ldots,u_N$ and $\theta$ such that the $u_i$ are mutually $(\ep,\theta)$-relatively close and $O_{u_i}=1$.
\end{lemma}
\begin{proof}
The statement is true for $N=2$ by the previous lemma. We will now prove that if the statement is true for $N$ then it is also true for $2N$.

So we now assume that the claim of the lemma is true for $N$. Find $u_1,\ldots,u_N$ and $\theta_1$ for $\ep_1 < \ep/6$.
Let $r_{u_{min}}=\min_{i=1,\ldots,N} r_{u_i}$. Choose $\ep_2>0$ so small that
\begin{equation}\label{ep2-1}
\ep_2 < \min \left\{e^{-\ep/6} \left( \frac{\ep}{3} - \ep_1\right) r_{u_{min}} , \frac{\ep}{6}\right\}
\end{equation}
and
\begin{equation}\label{ep2-2}
\left|1-e^{\ep_2}\right| + \ep_2 <\frac{\ep}{3} r_{u_{min}}.
\end{equation}

Consider the function $\phi(\theta) = \theta-\theta_1$. Apply Lemma~\ref{cl2of2} with $\ep_2$ and $\phi(\cdot)$ to find $s,t, \theta_2$ and $\widetilde{\om}$ satisfying conditions (i)-(v) of the Lemma, that is,
\begin{lista}
\item $\frac{r_s}{r_t} \in (e^{-\ep_2},e^{\ep_2})$;
\item $|\theta_s -\theta_t|<\ep_2$ and $O_s=O_t=1$;
\item $\Pi_{\theta_2}(s\widetilde{\om}) = \Pi_{\theta_2}(t\widetilde{\om})$;
\item $ |\theta_2 - (\theta_s + \theta_1) | < \ep_2$ and $ |\theta_2 - (\theta_t + \theta_1) | < \ep_2$.
\end{lista}
We now claim that the $2N$ distinct words $su_1,\ldots,su_N,tu_1,\ldots,tu_N$ are mutually $(\ep,\theta_2)$-relatively close. 

Since $\ep_1,\ep_2 < \ep/6$, it is easily seen that (i) and (ii) of Definition \ref{def_ep} are satisfied with $\ep/3$.  Clearly, $su_1,\ldots,su_N$ are mutually $(\ep_1,\theta_s+\theta_1)$-relatively close. We now prove that they are also $(\ep/3,\theta_2)$-relatively close:

Without loss of generality, we can assume $r_t\leq r_s$. Note $r_s <  e^{\ep/6}r_t$. Let $\om(u_i,u_j)$ be as described in part (iii) of Definition \ref{def_ep}. Denote by $\vec{e}_\theta$ the unit vector in the plane in the direction $\theta$. Observe that $|\vec{e}_{\al_1} -\vec{e}_{\al_2}| \leq |\al_1-\al_2|$. Then for any $i,j \in\{1,\ldots,m\}$ and $\om=\om(u_i,u_j)$, using (\ref{ep2-1}) and (d) we get (see Fig. \ref{resim1})

\begin{gather}
\begin{split}\label{long_comp}
|\Pi_{\theta_2}(s u_i & \om) - \Pi_{\theta_2}(s u_j\om)|\\
& =|(\Pi(s u_i \om) - \Pi(s u_j\om))\cdot \vec{e}_{\theta_2}|\\
&= |(\Pi(s u_i \om) - \Pi(s u_j\om))\cdot \vec{e}_{\theta_s + \theta_1} + (\Pi(s u_i \om) - \Pi(s u_j\om))\cdot (\vec{e}_{\theta_2} - \vec{e}_{\theta_s + \theta_1})|\\
& \leq |\Pi_{\theta_s+\theta_1}(s u_i \om) - \Pi_{\theta_s+\theta_1}(s u_j\om)| + |(\Pi(s u_i \om) - \Pi(s u_j\om))\cdot (\vec{e}_{\theta_2} - \vec{e}_{\theta_s + \theta_1})|\\
& \leq  |\Pi_{\theta_s+\theta_1}(s u_i \om) - \Pi_{\theta_s+\theta_1}(s u_j\om)| + |\Pi(s u_i \om) - \Pi(s u_j\om)| |\theta_2-(\theta_s + \theta_1)| \\
&  \leq \ep_1 D \min\{r_{su_i},r_{su_j}\} + D r_s \ep_2 \\
& < \ep_1 D \min\{r_{su_i},r_{su_j}\} + D e^{\ep/6}r_t  e^{-\ep/6} \left(\frac{\ep}{3} - \ep_1\right) r_{u_{min}} \\
& \leq \frac{\ep}{3} D \min\{r_{su_i},r_{su_j}\}.
\end{split}
\end{gather}

Therefore $su_1,\ldots,su_N$ are mutually $(\frac{\ep}{3},\theta_2)$-relatively close. A similar proof shows that $tu_1,\ldots,tu_N$ are also mutually $(\frac{\ep}{3},\theta_2)$-relatively close. Let $i_0$ be such that $r_{u_{i_0}}=r_{u_{min}}$. Now we will prove that $su_{i_0}$ and $tu_{i_0}$ are $(\frac{\ep}{3},\theta_2)$-relatively close and that we can use any $\om\in\Om$ in Definition \ref{def_ep} (iii). From this the result will follow, since for $\om=\om(u_i,u_j)$ we get

\[
\begin{split}
 \left|\Pi_{\theta_2}(s u_i \om) - \Pi_{\theta_2}(t u_j\om)\right|  \leq &
\left|\Pi_{\theta_2}(s u_i \om) - \Pi_{\theta_2}(s u_{i_0}\om)\right| 
+\left|\Pi_{\theta_2}(s u_{i_0} \om) - \Pi_{\theta_2}(t u_{i_0}\om)\right| \\
& + \left|\Pi_{\theta_2}(t u_{i_0} \om) - \Pi_{\theta_2}(t u_j\om)\right|
\end{split}
\]
(we can use any $\om\in\Om$ as $\om(u_i,u_j)$ if $i=j$).
\begin{figure}[h]
\psfrag{label1}{$l_{\theta_s + \theta_1}$}
\psfrag{label2}{$l_{\theta_2} = l_\theta$}
\psfrag{label3}{$\leq \ep_2$}
\psfrag{label4}{$\leq \ep_1 D r_s \min\{r_{u_i}, r_{u_j}\}$}
\psfrag{label5}{$\leq D r_s$}
\psfrag{label6}{$\Pi(s u_i \om)$}
\psfrag{label7}{$\Pi(s u_j \om)$}
\begin{center}\scalebox{.8}{\includegraphics{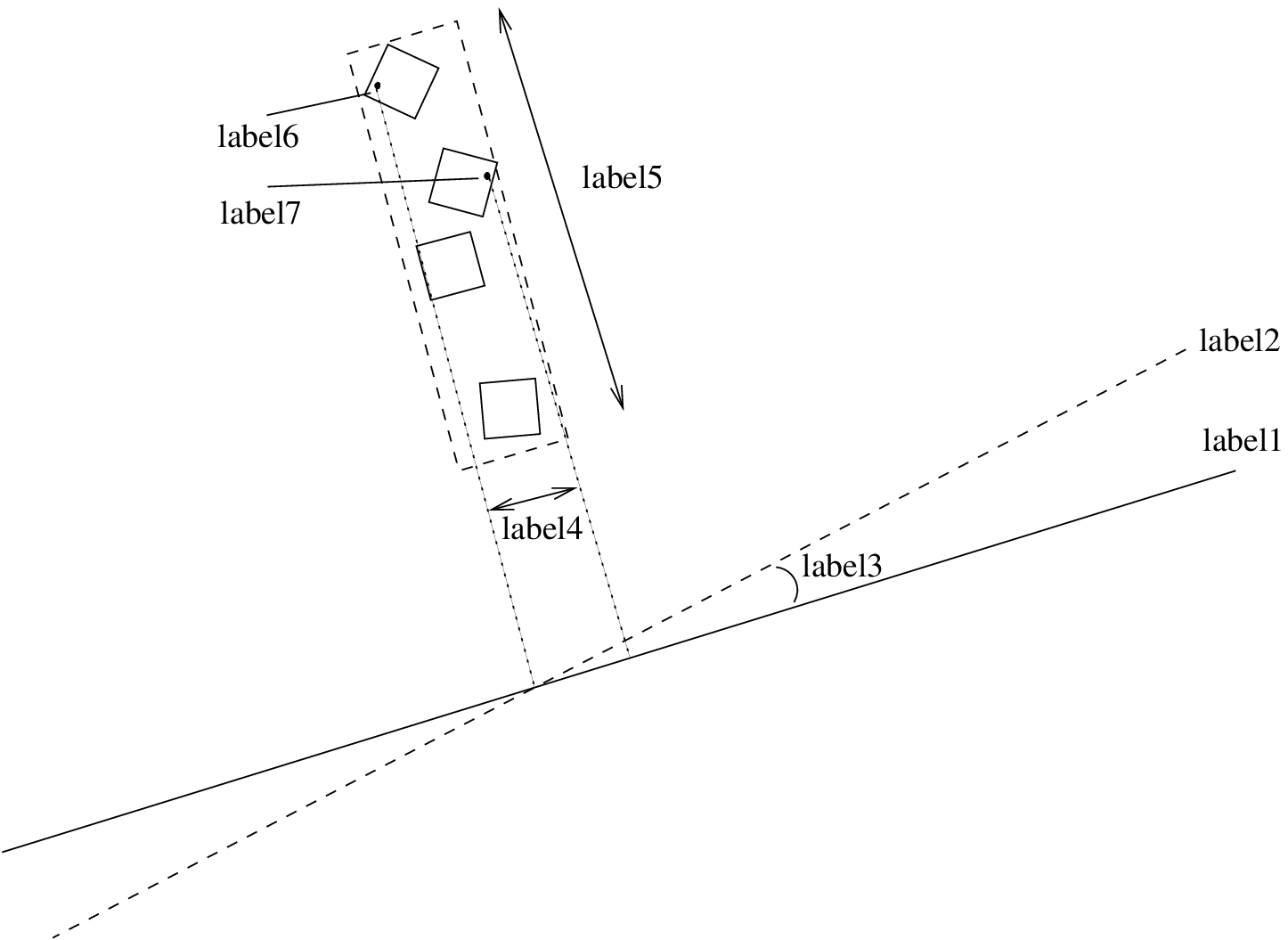}}\end{center}
\caption{}
\label{resim1}
\end{figure}

\begin{figure}[h]
\psfrag{label8}{$l_{\theta_2}$}
\psfrag{x}{$x^\prime$}
\psfrag{y}{$y^\prime$}
\begin{center}\scalebox{.75}{\includegraphics{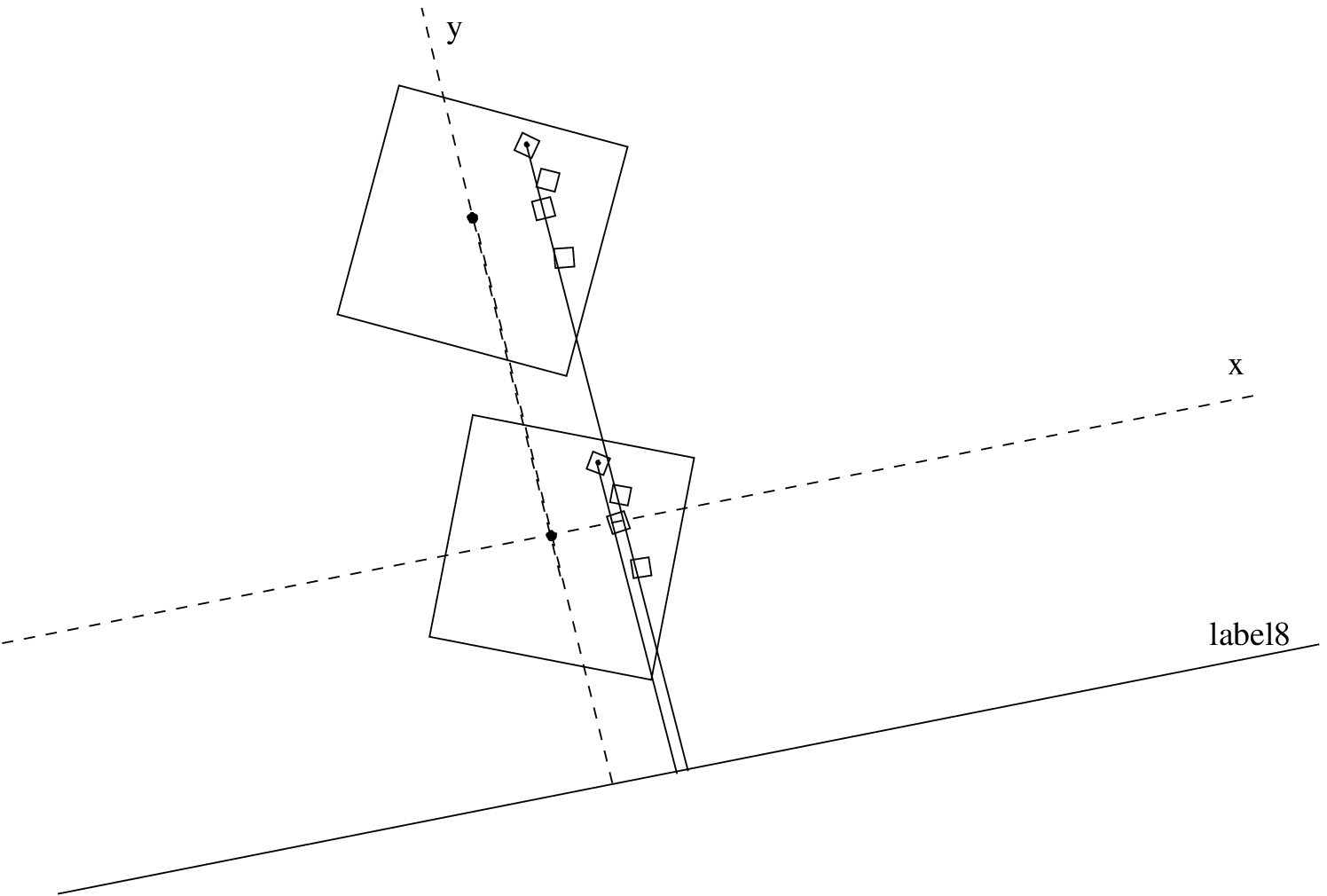}}\end{center}
\caption{}
\label{resim2}
\end{figure}

Denoting by $R_\theta$ the rotation map in the plane by angle $\theta$, for any word $u$ with $O_u=1$ we can write
\[
F_u (z) = r_u R_{\theta_u} z + b_u
\]
for some vector $b_u$. Also recall that $\Pi(u\tau)=F_u(\Pi(\tau))$ for any $\tau\in\Om$ by definition. Using (\ref{ep2-2}), the conditions (a)-(c) above and the linearity of $R_\theta$ we get (see Fig. \ref{resim2})
\begin{eqnarray*}
|\Pi_{\theta_2}(su_{i_0}\om)-\Pi_{\theta_2}(tu_{i_0}\om)|
& = & |\Pi_{\theta_2}(su_{i_0}\om)-\Pi_{\theta_2}(tu_{i_0}\om) + \Pi_{\theta_2}(t\widetilde{\om}) -\Pi_{\theta_2}(s\widetilde{\om})|\\
& = & |(\Pi(su_{i_0}\om)- \Pi(tu_{i_0}\om) + \Pi(t\widetilde{\om}) -\Pi(s\widetilde{\om})) \cdot \vec{e}_{\theta_2}| \\
& \leq & |\Pi(su_{i_0}\om)- \Pi(tu_{i_0}\om) + \Pi(t\widetilde{\om}) -\Pi(s\widetilde{\om})| \\
& = & |r_s R_{\theta_s}\Pi(u_{i_0}\om) - r_t R_{\theta_t}\Pi(u_{i_0}\om) \\
& & +\  r_t R_{\theta_t}\Pi(\widetilde{\om})- r_s R_{\theta_s}\Pi(\widetilde{\om}) |\\
& = & |r_s R_{\theta_s} (\Pi(u_{i_0}\om) - \Pi(\widetilde{\om})) - r_t R_{\theta_t} (\Pi(u_{i_0}\om) - \Pi(\widetilde{\om}))|\\
& \leq &  |(r_s - r_t) R_{\theta_s} (\Pi(u_{i_0}\om) - \Pi(\widetilde{\om}))|\\
& & + r_t |(R_{\theta_s}-R_{\theta_t}) (\Pi(u_{i_0}\om) - \Pi(\widetilde{\om}))|\\
& \leq & |r_s - r_t| \| R_{\theta_s}\| \cdot D + r_t \|R_{\theta_s}-R_{\theta_t}\|\cdot D\\
& \leq & D r_t (|1-e^{\ep_2}|+\ep_2) \leq D r_t \frac{\ep}{3} r_{u_{min}}\\
& \leq & \frac{\ep}{3} D r_{t u_{i_0}}
\end{eqnarray*}
and the result follows since $r_{t u_{i_0}} \leq \min\{r_{su_i},r_{tu_j}\}$ for any $i,j$. The lemma is proved.

\end{proof}

And finally we prove the following result which immediately implies the main result of this section:

\begin{proposition}\label{cl4of2}
For all $\theta$ we have $\Ha^\gamma (\Pi_\theta \Om) =0$.
\end{proposition}
\begin{proof}
Let $\mu_\theta$ be the projection of the measure $\mu$ under $\Pi_\theta$. We will prove that given any $N$, at $\mu_\theta$-almost all points $x$ the $\gamma$-dimensional upper density of $\mu_\theta$ given by
\[
\limsup_{r\to 0} \frac{\mu_\theta B(x,r)}{(2r)^\gamma}
\]
is at least $cN$ , where $c$ is a constant independent of $N$ (here $B(x,r)$ denotes an open ball of radius $r$ around $x$) . The result will follow from a standard density theorem (see \cite{Fal2}~Prop.~2.2  for a statement).

Observe that if $u,v$ are distinct and $(\ep,\theta)$-relatively close for some $\ep$ and $\theta$, then $u$ and $v$ can not be subwords of each other provided $\ep$ is small enough. This follows from the fact that $r_u / r_v \notin (r_{min},1/r_{min})$ if any of $u$ and $v$ is a proper subword of the other. Note that this also implies that $[u]$ and $[v]$ are disjoint in $\Om$. Now, given $N$, we can find distinct $u_1,\ldots,u_N$ that are mutually $(1,\theta_0)$-relatively close for some $\theta_0$. This we can see by applying  Lemma~\ref{cl3of2} with any $\ep<1$. By using small enough $\ep$ in the Lemma if necessary, we can assume that $[u_i]$ are disjoint. For the purposes of our proof, the words $u_i$ will be regarded as mutually $(1,\theta_0)$-relatively close.

Now let $x=\Pi_\theta (\tau)$ where $\tau$ satisfies the conditions of Lemma~\ref{cl1of2} with $\ep= r_{u_1}$, word $u_1$ and $\phi=\theta - \theta_0$. That is, $\tau$ has prefixes of the form $s_j u_1$ where $s_j$ satisfy $|\theta_{s_j} + \theta_0 - \theta| < r_{u_1}$. By Lemma~\ref{cl1of2}, $\mu_\theta$-almost all $x$ are of this form. It follows that for each $j$,$i$ and $\om\in\Om$, using ideas similar to those in (\ref{long_comp}) we get

\begin{eqnarray*}
| x - \Pi_\theta (s_j  u_i \om) | &= & |  \Pi_\theta (\tau) - \Pi_\theta (s_j u_i \om)| \\
 &= & | (\Pi (\tau) - \Pi (s_j u_i \om)) \cdot \vec{e}_\theta|\\
&\leq & | (\Pi (\tau) - \Pi (s_j u_i \om)) \cdot \vec{e}_{(\theta_{s_j} + \theta_0)}| \\
& &+ \ | (\Pi (\tau) - \Pi (s_j u_i \om)) \cdot (\vec{e}_{\theta} - \vec{e}_{(\theta_{s_j} + \theta_0)})|\\
&\leq & |\Pi (\tau) - \Pi (s_j u_i \om)| + | \Pi (\tau) - \Pi (s_j u_i \om)| |\theta - (\theta_{s_j} + \theta_0)|\\
&\leq & |\Pi (\tau) - \Pi (s_j u_1 \om)| + | \Pi (s_j u_1 \om)- \Pi (s_j u_i \om)| \\
&&+\  D r_{s_j}|\theta - (\theta_{s_j} + \theta_0)|\\
&\leq & D r_{s_j u_1} + D r_{s_j} \min\{r_{u_1},r_{u_j}\} + D r_{s_j} r_{u_1}\\
&\leq & 2 D r_{s_j u_1} + D r_{s_j} r_{u_1} = 3 D r_{s_j u_1}.  
\end{eqnarray*}

\begin{figure}[h]
\psfrag{label10}{$l_{\theta_{s_j} + \theta_0}$}
\psfrag{label9}{$l_\theta$}
\psfrag{label11}{$\leq D r_{s_j}$}
\psfrag{label12}{$\leq D r_{s_j} \min_i \{r_{u_1}, r_{u_i}\} + D r_{s_j} r_{u_1} \leq 2 D r_{s_j u_1}$}
\psfrag{label13}{$x$}
\psfrag{label14}{$s_j$}
\begin{center}\scalebox{.8}{\includegraphics{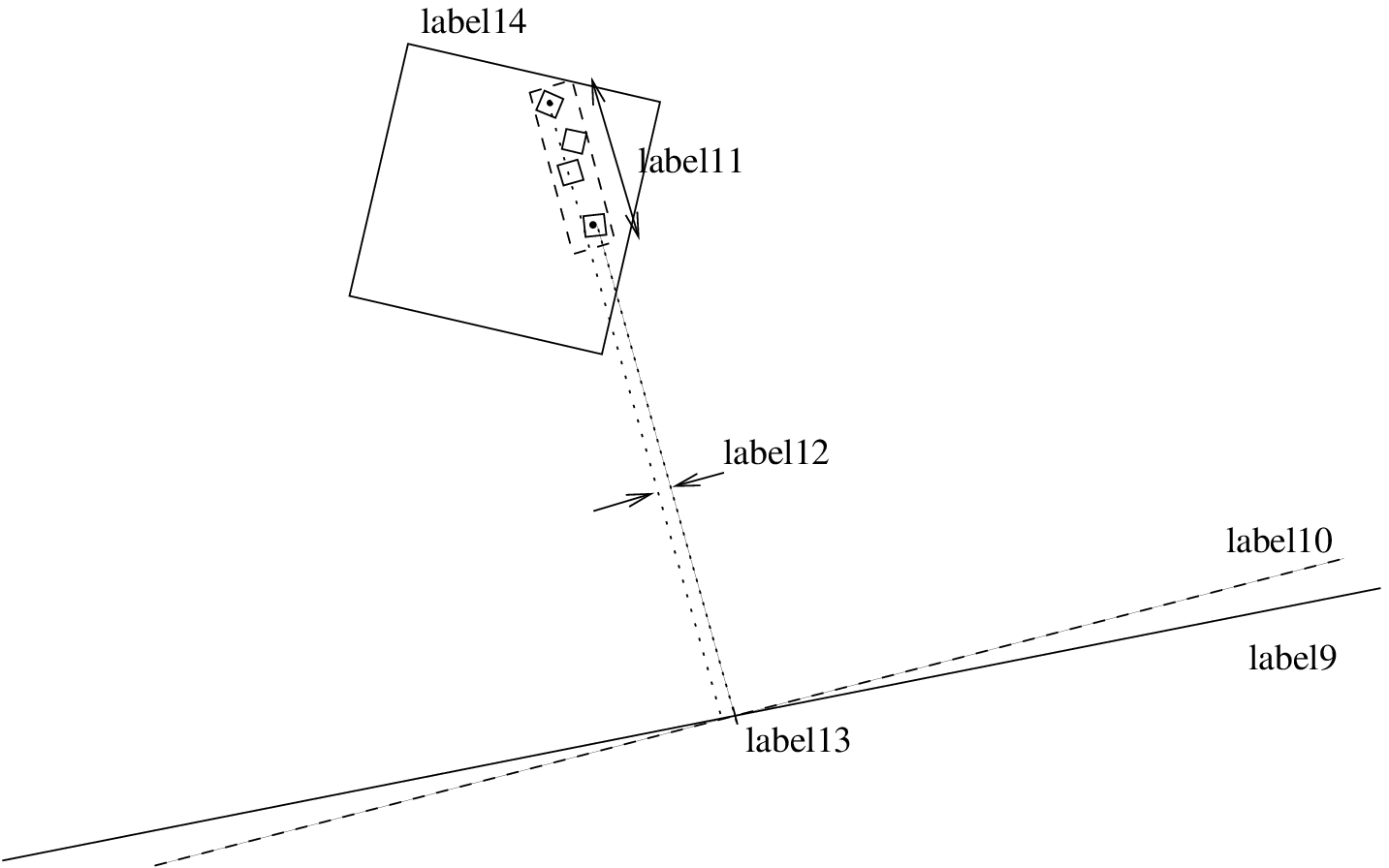}}\end{center}
\caption{}
\label{resim3}
\end{figure}

Thus if we set $b_j =5 D r_{s_j u_1}$ then each $\Pi_{\theta}(s_j u_i)$ lies in $B(x,b_j)$, thus (using that $r_{s_j u_i} \geq \frac{1}{e} r_{s_j u_1}$ for each $i$):

\[
\frac{\mu_\theta (B(x,b_j))}{b_j^\gamma} \geq \frac{1}{(10De)^\gamma} N
\]

and this proves our claim.
\end{proof}

\bigskip
\noindent {\bf Applications to visibility.} We now make remarks about the application of the method above to visibility problems. We consider the visibility of a self-similar set from a point. Given a point $a\in\mathbb{R}$, we define the radial projection (centered at $a$) as
\[
\pi_a : \mathbb{R}^2\setminus {a} \longrightarrow S^1, \quad \pi_a(x)=\frac{x-a}{|x-a|}.
\]
\begin{definition} 
We say that a set $E$ is $s$-visible from $a$ if $\Ha^s(\pi_a(E))>0$. We say $E$ is $s$-invisible from $a$ if $\Ha^s(\pi_a(E))=0$.
\end{definition}
\noindent \textit{Remark.} The standard terminology is to say ``visible/invisible from $a$'' in the case $s=1$.

We observe that if $E$ satisfies the conditions of Theorem \ref{Main2}, we can use the radial projections of the measure $\mu$ to show that, at all points, the projected measure has infinite $\gamma$-dimensional density almost everywhere (where $\gamma$ is the similarity dimension). More precisely, given $a\in\mathbb{R}^2$, we define $\tilde{\mu}_a$ to be the restriction of $\mu$ to $\Om\setminus \Pi^{-1}(a)$ and define $\mu_a$ as the projection of $\tilde{\mu}_a$ to a measure on $S^1$ via $\pi_a \circ \Pi$. It is easy to see that, in the proof of Lemma~\ref{cl4of2}, one can use Corollary~\ref{corol_to_1of2} and the same idea of aligning the 1-relatively close squares around typical points in such a way that the radial projections of these squares have large overlap, making the measure density big. Since the required modifications are fairly obvious, we state the result without a proof:

\begin{theorem}\label{Main2v}
Let $E\subset\mathbb{R}^2$ be a self-similar set for which $\{\theta_u:\text{$u$ is a word}\}$ is dense modulo $2\pi$. If $\gamma$ is the similarity dimension of $E$, then for all points $a$ in the plane, $E$ is $\gamma$-invisible from $a$.
\end{theorem}

We also remark that in a recent work \cite{Sol_Karol} of K.~Simon and B.~Solomyak, it is proven that purely unrectifiable planar self-similar sets with finite $\mathcal{H}^1$ measure and satisfying the open set condition are invisible from all points in the plane. This implies our conclusion in the case $d=1$. To see this, note that if $d=1$ and $\mathcal{H}^1 (E)>0$, dense set of rotations implies the nonexistence of tangent directions at all points of $E$, hence the pure unrectifiability of $E$ (see \cite{Fal1} for an overview of tangency properties).

\section{Favard length in some special cases}\label{sec_favard}

The Favard length of a planar set $E$ is defined as
\[
\fav E = \int_0^\pi \mathcal{L}^1 (E^{\theta}) d\theta
\]
where $E^\theta$ is the orthogonal projection of $E$ onto $l_\theta$. It can be interpreted as a measure of the probability of Buffon's needle hitting the set $E$. By Besicovitch's well-known theorem, irregular 1-sets project to zero measure in Lebesgue almost all directions, therefore their Favard length is zero (see \cite{Fal1} for details). One can then ask about the speed at which $\fav E(\rho)$ decreases to zero, where $E(\rho)$ denotes the $\rho$-neighborhood of the set $E$.

Mattila has shown \cite{Mattila_Favard} that $c/(-\log \rho)$ is a lower bound for 1-sets. The question about the best possible (general) upper bound is open. We should remark that different upper bounds may apply to different classes of sets.

In \cite{Sol_Buffon}, Y.~Peres and B.~Solomyak gave an upper bound for self-similar 1-sets in the plane that do not involve rotations and satisfy the strong separation condition. Defining $\log_* x$ as in (\ref{def_log_star}), they proved that for such sets $C \exp(-a \log_* (1/\rho))$ for some $a, C>0$ is an upper bound. This bound is of course far from the lower bound given by Mattila, and it is not known whether it can be improved or not. Peres and Solomyak also gave an example of a random Cantor set for which the expected upper bound is of the same order as Mattila's lower bound. It is an interesting problem to give more accurate estimates for deterministic sets. The method used in this section is based on the approach in \cite{Sol_Buffon}.

Now we return to our second main result. We are going to consider a homogeneous self-similar set $E$ in the plane of similarity dimension $1$, defined by $m$ maps. Let $r=1/m$ be the common contraction factor of the homogeneous system. We are assuming that the defining maps produce two non-rotating maps of the same contraction rate and a map containing rotation by a Diophantine multiple of $2\pi$. By composing the non-rotating maps with themselves to remove reflections if necessary, we can assume that
\begin{lista}
\item There are two distinct words $u,v$ with $|u|=|v|=k$, $O_u = O_v =1$ and $\theta_u = \theta_v=0$;
\item One of the defining maps, say, $F_1$, contains irrational rotation by $\theta_1$, where $\theta_1 / 2\pi$ is $(c,d)$-Diophantine.
\end{lista}

Under these assumptions, we are going to prove the following:

\begin{theorem}\label{MainFav}

Given any $\delta>0$, there exists $A>0$ such that
\begin{equation}\label{eqn_favard}
\mathcal{L}^1 (E^{\theta}(\rho)) \leq \frac{A}{(\log(-\log \rho))^B}
\end{equation}
uniformly for all $\theta$, where $B=\log 2 / ((1+\delta)k(d+1) \log m)$. Thus,
\[
\fav (E(\rho)) \leq \frac{\pi A}{(\log(-\log \rho))^B}.
\] 
\end{theorem}

Observe that it suffices to prove this theorem for a sequence of $\rho_n $ decreasing to $0$ such that $\log(-\log \rho_{n+1}) / \log(-\log \rho_n)$ is bounded. Also for simplicity, we will assume that the diameter of $E$ is $1$ (this will only change the constant $A$). For a word $u$, we have $r_u=r^{|u|}$ since the system is homogeneous.

Before we begin the proof, we mention an equivalent reformulation of this result: Let $C$ be the convex closure of the self-similar set $E$. Let

\[
C_n = \bigcup_{u:\ |u|=n} F_u (C).
\]
Observe that there exists a constant $K$ independent of $n$ and $K$ translation maps in the plane such that, for $\rho\leq r^n$, $E(\rho)$ can be covered by these $K$ translates of $C_n$. This follows from the fact that $C$ contains a (nontrivial) ball since there are irrational rotations. Therefore, $\fav(C_n)$ is comparable (with uniform constants) to $\fav(E(r^n))$. Note that if we take $\rho = r^n$, we have $n=\log r / (-\log \rho)$ therefore (\ref{eqn_favard}) becomes

\begin{equation}\label{eqn_favard2}
\mathcal{L}^1 (C_n^\theta) \leq \frac{\tilde{A}}{(\log n)^B}.
\end{equation}
In this formulation, the lower bound given for $\fav(C_n)$ by Mattila's result is $c/n$.

\medskip
The first stage in the proof is to construct $2^n$ $(0, \theta)$-relatively close words for any given $n$: Let $u,v$ as in condition (a) given above. Let $k=|u|=|v|$. Note that $u$ and $v$ satisfy the definition of $(0,\theta)$ relative closeness if we choose $\om= \om(u,v)= \bar{u}$ and $\theta$ to be perpendicular to the line connecting $\Pi(u\om)$ to $\Pi(v\om)$ (or any line if these points coincide). Now we observe that $uu, uv,vu,vv$ are also mutually $(0,\theta)$ relatively close words with the same $\theta$ and $\om$ used for $u$ and $v$: Clearly $\Pi(u\om) = \Pi(uu\om)$ and $\Pi(v\om) = \Pi(vu\om)$. And since $F_u, F_v$ contain no rotation, the points $\Pi(uv\om)$ and $\Pi(vv\om)$ lie on the line connecting $\Pi(uu\om)$ to $\Pi(vu\om)$. Continuing this procedure, for any $n$, we can obtain $2^n$ words $u_1,\ldots,u_{2^n}$ that are mutually $(0,\theta)$-relatively close and with $|u_i|=kn$ for all $i=1,\ldots,2^n$.

By our assumption $\theta^\prime := \theta_1 / 2\pi$ is $(c,d)$-Diophantine, that, is for any integers $N,M$ we have
\begin{equation}\label{eqn_dioph}
|N \theta^\prime - M| > c N^{-d}.
\end{equation}
Observe that, given any $\ep>0$, by the pigeonhole principle, there are integers $N\leq \frac{1}{\ep}$ and $M$ such that 
\[
|N \theta^\prime - M| < \ep.
\]
Combined with (\ref{eqn_dioph}), we have
\begin{equation}\label{eqn_dioph2}
\ep> |N \theta^\prime - M| > c \ep^d.
\end{equation}
Now we make the following observations: 

{\bf Observation 1.} There exists a constant $c_1=c_1(\theta^\prime)$ such that, given any $\ep$, there is an integer $p< c_1 \ep^{-d-1}$ such that the numbers

\[
\theta_1, \theta_{11}, \theta_{111},\ldots,\theta_{1^p}
\]
form an $\ep$-net modulo $2\pi$. Therefore, given any word $s$, angle $\phi$ and number $\ep$, there is a word $t_s=1^p$ with $|t_s| < c_1 \ep^{-d-1}$ such that $|\theta_{s t_s} - \phi|< \ep$.

{\bf Observation 2.} Simple geometric arguments similar to those in the previous section (e.g.~ see the computation in (\ref{long_comp})) show that, given $(0,\theta)$-relatively close cylinders $u_1,\ldots,u_{2^n}$, if the word $s$ is such that $O_s=1$ and \mbox{$|\theta_s + \theta_{u_1} - \theta_0| < r^{|u_1|}$}, then the cylinders $su_1,\ldots,su_{2^n}$ are $(1,\theta_0)$-relatively close.

Now we fix an $n$ and the corresponding $u_i$, with $|u_i|=kn$. Let 
\begin{equation}\label{eqn_sn}
\begin{split}
&s(n)= 2 c_1 m^{(d+1)kn} \\
&L_n = m^{s(n)} s(n)^2\\
&\rho_n = r^{L_n}.
\end{split} 
\end{equation}

Given $\theta_0$, our purpose is to prove (\ref{eqn_favard}) with $\theta=\theta_0$, $\rho=\rho_n$, and $n$ sufficiently large. Note that $\log(-\log \rho_{n+1}) / \log(-\log \rho_n)$ is bounded, therefore this suffices for the general proof.

Let $G$ be the set of words of length $L_n$. We partition $G$ into two subsets: Those words that contain $u_1$ ``at the right place'' in the sense of Lemma \ref{cl1of2}, and others. Let

\[
 G_1 = \{\al\in G\mid \exists s:\ \text{$s u_1$ is a prefix of $\al$ with } O_s=1, |\theta_s + \theta_{u_1} - \theta_0|< r^{|u_1|}\}
\]
and define $G_2=G\setminus G_1$. For a word $\al$ denote $E_\al = F_\al (E)$. As before, we denote by $E_\al^{\theta_0} (\rho)$ the projection of $E_\al(\rho)$ onto $l_{\theta_0}$. Then

\begin{equation}\label{eqn_h1_h2}
E^{\theta_0}(\rho_n) = \bigcup_{\al \in G_1} E_\al^{\theta_0} (\rho_n)  \ \cup\ \bigcup_{\al \in G_2} E_\al^{\theta_0} (\rho_n) =: H_1 \cup H_2.
\end{equation}

We first give an upper bound for the number of words in $G_2$:

By the observations above, there is a word $t$ of length no more than $c_1 (r^{|u_1|})^{-d-1} = s(n)/2$ such that $O_t=1$ and $|\theta_t + \theta_{u_1} - \theta_0|< r^{|u_1|}$. Then we have $|tu_1|\leq s(n)/2 + kn \leq s(n)$, for $n$ sufficiently large. Write $\Om \setminus [tu_1]$ as a disjoint union of cylinders represented by words of length $|tu_1|$. Now for each cylinder $[s]$ in this union, we can find a suffix $t_s$ of length no more than $s(n)/2$ that satisfies $O_{s t_s}=1$ and $|\theta_{s t_s} + \theta_{u_1} - \theta_0|< r^{|u_1|}$. Then $|t_s u_1| \leq s(n)$, and we repeat this procedure by writing $[s]\setminus [s t_s u_1]$ as a disjoint union of cylinders represented by words of length $|s t_s u_1| \leq |s|+s(n)$ for each $s$. Note that at each stage we remove words that are in $G_1$. After $L_n / s(n)$ steps we will have obtained a collection of words of length no more than $L_n$ that contain all possible prefixes for the words in $G_2$. Since at each step we disallow non-overlapping subwords of length $\leq s(n)$, we have

\begin{equation}\label{eqn_exps}
\#G_2 \leq (m^{s(n)}-1)^{L_n/s(n)} = m^{L_n} \left( 1 - \frac{1}{m^{s(n)}} \right) ^ {m^{s(n)} s(n)} \leq m^{L_n} e^{-s(n)}.
\end{equation}

Recall that $\diam E_\al = r^{|\al|}$ with $r=\frac{1}{m}$. Then, in the view of (\ref{eqn_sn}), we get

\begin{equation}\label{eqn_h2}
\mathcal{L}^1(H_2) \leq \sum_{\al \in G_2} \diam E_\al^{\theta_0} (\rho_n) \leq (\# G_2) 3 r^{L_n} \leq 3 e^{-s(n)}.
\end{equation}

Now we turn our attention to the words in $G_1$. A word in $G_1$ is of the form $su_1\beta$ where $O_s=1$ and $|\theta_s + \theta_{u_1} - \theta_0|< r^{|u_1|}$. Then the words $su_1,\ldots,su_{2^n}$ are mutually $(1,\theta)$-relatively close by Observation 2. Also recall that for any word $\al$, $\mu([\al]) = r^{|\al|} = \diam E_\al$. Therefore, we can find a constant $\eta$ (not depending on $n$ or $x$) such that for each $x\in H_1$, there is a ball $B(x,R_x)$ with

\[
\mu_{\theta_0} (B(x,R_x)) \geq 2^n \eta R_x.
\]
This can be proved by arguing as in the proof of Proposition \ref{cl4of2}. Indeed, if $x\in \Pi([su_1])$ö we can take $R_x= 3 r^{|su_1|}$ and $\eta=1$. Then by Vitali Covering Theorem this implies that
\begin{equation}\label{eqn_h1}
\mathcal{L}^1(H_1) \leq \zeta 2^{-n}
\end{equation}
for some $\zeta$ not depending on $n$. Therefore, by (\ref{eqn_h1_h2}), our theorem will be proved if we can put upper bounds to $2^{-n}$ and $e^{-s(n)}$ in terms of $\rho=\rho_n$ as in (\ref{eqn_favard}).

By (\ref{eqn_sn}), we have

\[
\log \rho_n = -L_n \log m = -m^{s(n)} s(n)^2 \log m
\]
thus
\[
\log (- \log \rho_n ) =  s(n)\log m + 2 \log s(n) + \log\log m \leq(1+\log m) s(n)
\]
for $n$ sufficiently large, which implies
\begin{equation}\label{eqn_sn_bigger}
s(n) = 2c_1 m^{kn(d+1)} \geq \frac{1}{1 + \log m} \log (-\log \rho_n).
\end{equation}
Taking logs above, we get
\[
\log(2c_1) + k(d+1) n \log m \geq \log(1/(1 + \log m)) + \log\log (-\log \rho_n).
\]
which implies that, for any $\delta>0$
\begin{equation}\label{eqn_n_bigger}
n \geq \frac{1}{(1+\delta)k(d+1)\log m} \log\log (-\log \rho_n).
\end{equation}
for $n\geq N(\delta)$. Setting $B= \log 2 / ((1+\delta)k(d+1) \log m)$ we get
\[
2^{-n} =\frac{1}{e^{n \log2}} \leq \frac{1}{e^{B \log\log (-\log \rho_n).}} = \frac{1}{(\log (-\log \rho_n))^B}.
\]
And finally, by (\ref{eqn_sn_bigger}), we have (for large $n$)
\[
e^{-s(n)} \leq \frac{1}{e^{\frac{\log (-\log \rho_n)}{1+\log m}}} = \frac{1}{(-\log \rho_n)^{1/(1 + \log m)}} \leq \frac{1}{(\log (-\log \rho_n))^B}.
\]
The proof is complete.

\bigskip
\bigskip

\noindent {\bf Acknowledgements.} The author would like to thank Boris Solomyak for his valuable comments and suggestions.

\bibliographystyle{hamsplain}
\bibliography{general}

\end{document}